\documentclass[11pt]{article}
\usepackage{times}
\usepackage[english]{babel}
\usepackage{subfigure}
\usepackage{epsfig,graphicx,psfrag}
\usepackage{amsfonts,amsmath,amssymb, amsthm}
\usepackage{fullpage}						
\usepackage{hyperref}

\newtheorem{theorem}{Theorem}

\newtheorem{remark}[theorem]{Remark}

\def\trds{\,\textrm{ds}}

\def\trdt{\,\textrm{dt}}
\def\trdBt{\,\textrm{dB}_{{\rm t}}}

\newcommand*{\argmin}{\operatornamewithlimits{argmin}\limits}
\newcommand*{\argmax}{\operatornamewithlimits{argmax}\limits}

\begin{document}
\title{\huge Stochastic Optimal Control for Online Seller under Reputational Mechanisms}
\author{
	Milan Bradonji\'c\thanks{Milan Bradonji\'c is with Mathematics of Networks and Systems, Bell Labs, Alcatel-Lucent, 600 Mountain Avenue, Murray Hill, NJ 07974, USA, e-mail: \texttt{milan@research.bell-labs.com}.}, 
	Matthew Causley\thanks{Matthew Causley is with the Department of Mathematics, Kettering University, Flint, MI 48504, USA, e-mail: \texttt{mcausley@kettering.edu}.}, 
	Albert Cohen\thanks{Albert Cohen is with the Department of Statistics and Probability, and Department of Mathematics, Michigan State University, East Lansing, MI 48824, USA, e-mail: \texttt{albert@math.msu.edu}.}
}

\maketitle

\begin{abstract}
In this work we propose and analyze a model which addresses the pulsing behavior of sellers in an online auction (store). This pulsing behavior is observed when sellers switch between advertising and processing states. We assert that a seller switches her state in order to maximize her profit, and further that this switch can be identified through the seller's reputation. We show that for each seller there is an optimal reputation, \emph{i.e.}, the reputation at which the seller should switch her state in order to maximize her total profit.
 We design a stochastic behavioral model for an online seller, which incorporates the dynamics of resource allocation and reputation. The design of the model is optimized by using a stochastic advertising model from~\cite{sethi-1983-deterministic} and used effectively in the Stochastic Optimal Control of Advertising~\cite{raman-2006-boundary}. This model of reputation is combined with the effect of online reputation on sales price empirically verified in~\cite{mink-2006-reputation}. We derive the Hamilton-Jacobi-Bellman (HJB) differential equation, whose solution relates optimal wealth level to a seller's reputation. We formulate both a full model, as well as a reduced model with fewer parameters, both of which have the same qualitative description of the optimal seller behavior. Coincidentally, the reduced model has a closed form analytical solution that we construct. 
\end{abstract}


\section{Introduction}
\label{sec:intro}

In this work, we consider sellers in an online market like Amazon without a bidding mechanism. The setting may include concurrent online auctions on a site like eBay, for a product or service with only a ``buy-it-now'' option, under the assumption that buyers can compare reputation among online platforms. In this online environment, the buyer and seller have a certain amount of anonymity, which has implications on the fairness of the auction~\cite{beam98auctions}. One of the only counter measures to ensure fairness is the use of feedback forums, such as Amazon, eBay, and Yahoo!, at online auction sites. The feedback provided by customers provides a ranking system for the sellers, which in turn rewards fair sellers and penalizes unfair or unreliable sellers \cite{bajari04economic, resnick06experimental, ResnickZeckhauser2002}. Since it is in the sellers best interest to maximize their profits, it follows that a seller will seek to use reputation as a means to optimize their decision making ~\cite{HouserWooders2006, dellarocas05reputation}; e.g., for the optimal bidding strategies in sequential auctions, see~\cite{apt09optimal}.

Online auction sellers can be broadly categorized by two types of behavior:  those focused on advertising and customer service, which includes following up with past customers; and those who instead focus only on processing existing orders. While the latter behavior will likely lower the seller reputation on average, this immediate increase in wealth may offset the long-term damage from a lower reputation. With this in mind, it is compelling to ask whether there exists an optimal long term strategy, in which a seller attempts to maximize her profit by achieving a fixed reputation.

We propose a simple behavioral model to study the relationship between a seller's wealth and reputation in the auctions environment with only the ``buy-it-now'' option. Following the standard Nerlove-Arrow construction~\cite{nerlove-1962-optimal}, which was extended to a stochastic setting by Sethi~\cite{sethi-1983-deterministic}, the wealth and growth of the seller are described by continuous stochastic processes. 
Note that our model is a continuous approximation of a discrete auctioning process. The discrete model of the process may be less tractable but still very interesting and useful. 
We propose an optimal infinite horizon strategy, and incorporate the wealth-growth model of Mink and Seifert~\cite{mink-2006-reputation}, based on empirical studies of ebay sellers.  Also, we use the model for reputation proposed by Sethi~\cite{sethi-1983-deterministic} and Rao~\cite{rao-1986-estimating} which has been adopted by many others in optimal advertising models, such as Raman's recent work in Boundary Value Problems in the Stochastic Optimal Control of Advertising~\cite{raman-2006-boundary}. 
This work incorporates the idea that humans do not multi task well, but rather switch, see~\cite{book-levitin-2014}. In other words, switching costs do not factor in for an individual agent, but only for large online retailers who have marketing department reorganization costs. Additionally, the model in this work allows all sellers to have a little extra capacity/money for either shipping or advertising (e.g., writing a letter or similar individual campaign), which is rather small and not unlimited.

The rest of the paper is laid out as follows. In Section \ref{sec:mink.seifert}, we define the stochastic Nerlove-Arrow model, along with the growth rate due to Mink and Seifert. Using the principal of optimality we derive the corresponding Hamilton Jacobi Bellman (HJB) boundary value problem satisfied by the value function. In Section \ref{sec:mod.market.share}, we propose several important modifications to the existing methodology, and in doing so derive a new model which is normalized with respect to the seller reputation, and has a linear growth rate. In making these modifications, we also find a reduced version of our new model, which is governed by a geometric Brownian motion. We prove that both models admit unique piecewise continuous solutions, and are qualitatively similar. However the reduced model, surprisingly, has a closed form piecewise defined solution. This allows us to prove that pulsing behavior is an optimal strategy for sellers and also serves as a benchmark to our numerical solution of the full model in Section~\ref{sec:num.results}, which does not have a closed form analytic solution. Our results are compelling, and warrant further investigation for the finite horizon case, which we discuss with concluding remarks in Section \ref{sec:conclusion}.

\section{Explicit Resource Allocation Mechanism}
\label{sec:mink.seifert}

Consider a seller in an online marketplace as a triple $(W,R,\mu)$ representing her wealth, reputation, and excess rate, respectively. Here, the excess rate is extra effort that can be allocated to either expedite payment {\em or} increase her reputation. In this setting, the reputation $R$ is a positive number reflecting customer satisfaction. As an initial approach, we consider only the effects of reputation on the optimal way to grow her wealth, and leave the more general model with shipping costs for future work. Also, in our mathematical formulation,  we consider the reputation mechanism first suggested by the work of Nerlove and Arrow~\cite{nerlove-1962-optimal} and later generalized to stochastic settings in~\cite{sethi-1983-deterministic, rao-1986-estimating, raman-2006-boundary}.

The seller is one of many competing in the online marketplace, and has her transactions verified and processed automatically via the marketplace. All sellers competing in this marketplace have an expected speed for mailing out completed orders to be able to sell in the marketplace, with some time allowed for delayed shipping before the marketplace intervenes. In return, the marketplace guarantees the transaction. Payment is modeled to be released to the seller once the purchased item is mailed to the~buyer.

However, the seller can choose to go beyond these expected standards. For example, she may prepare packages and drive to the post office earlier than expected to expedite the payment from the marketplace. Or she may spend time and money communicating with buyers in a bid to increase their favorable ratings. This may include sending a small unexpected gift in the package to be mailed out or spending extra time to individually craft and send emails out to purchasers. Our model assumes that at most only one of these extra actions can happen with noticeable return at any given time. The seller may also decide to do nothing extra, and in this case $\mu = 0$.

For example, if she is focused on getting her payment released, then she is not able to focus on creative ways to engage her buyers post sale. We also assume that her extra capacity to do either extra action is limited by absolute minimum mailing rates and her pool of resources she is willing to contribute to going beyond expected standards.  By choosing to shift up her resources from promotion to expedited mailing and back, the seller can influence her wealth and reputation levels via the excess rate $\mu$. Positive $\mu$ corresponds to a promotional state, while negative $\mu$ corresponds to an expedited mailing state.

We note that if our seller is in fact a large company with many such agents that can work on either processing or sales, then the company can choose to shift the ratio of workers who work in processing or sales, and this model is the subject of future work. The work we present here is limited to a one agent seller. 
We note here that our model also assumes no cost for the individual to switch behavior from promotion to expedited mailing. In a large company, one would expect such costs to arise in restricting the size of advertising and shipping departments. One related paper that addresses such switching costs, albeit in a real options investment setting, is the paper by Duckworth \emph {et~al}.~\cite{duckworth-2001-model}.

We utilize a continuous time formulation for wealth and reputation, as well as control of these processes. The seller will be determined to behave optimally if she maximizes her expected present value of total earnings, until her reputational value reaches $R=0$. This produces a fully consistent formulation of both the value function $V$, and the strategy to achieve this optimal value.

\subsection{Mathematical Formulation}

Following the framework in \cite{book-fleming-soner-2006, fleming-1975-determinstic}, reputation $R$ is observed in $  \mathcal{O} := (0,\infty)$ and the control $\mu $ is constrained to $U :=[- \epsilon, \epsilon]$. Since this process 
will be stochastic, we introduce a probability space $(\Omega,\mathcal{F}, \mathbb{P})$ and corresponding reference 
probability systems $\nu = \left(\left(\Omega,\left \{ \mathcal{F}_{s} \right\}, \mathbb{P}\right), B \right)$ where $\mathcal{F}_{s} \subset \mathcal{F}$ and $B$ is an $\mathcal{F}_{s}$-adapted Brownian motion. The evolution of wealth and reputation is hence modeled  by a controlled two-state Markov process $(W^{\mu_{.}}_{.},R^{\mu_{.}}_{.})$, where wealth $W$ grows at a rate proportional to a function $ h(R)$ that links reputation to revenue per sale, subject to a controlled extra processing rate $\mu$  per unit time chosen from the admissible class $A_{\nu}$,

\begin{equation}
\mathcal{A}_{\nu} :=  \left \{ \mu \in \mathcal{F}_{s} \mid \mu \text{ is } \mathcal{F}_{s} \text{-p.m.}, \mu_{.} \in U \text{ on } [0,\infty) , \mathbb{E} \Big[ \int_{0}^{\tau_{\mathcal{O}}} e^{-\rho s} \left| (1-\mu_{s}) h(R^{(\mu)}_{s}) \right | \trds  \Big] < \infty \right \} \,. 
\end{equation}

Here, \text{p.m.} denotes {\em progressively measurable}, $\rho$ is the constant discount factor to account for the time value of money, and $\tau_{\mathcal{O}}$ is the first exit time of reputation $R$ from $\mathcal{O}$, or $\infty$ if $R \in \mathcal{O}$ for all $s \geq 0$. We set $\infty$ as an absorbing state for reputation. Symbolically, 

\begin{equation}
 \tau_{\mathcal{O}} := \inf \left \{ t \geq 0 \mid R^{(\mu)}_{t} \notin \mathcal{O} \right \} \,.
\end{equation}

Formally, our two state controlled Markov process is
\begin{eqnarray}
\nonumber
dW^{(\mu)}_{t} &=& (1-\mu_{t}) h(R^{(\mu)}_{t}) \trdt \\
\nonumber
dR^{(\mu)}_{t} &=& (\mu_{t} - \kappa R^{(\mu)}_{t}) \trdt + \sigma \trdBt \,,
\end{eqnarray}
where $\kappa$ is a proportionality constant which accounts for mean reversion~\cite{nerlove-1962-optimal} and $\sigma$ is the constant volatilty accounting for random effects to reputation.

The major difficulty is in defining the growth rate $h(R)$. However, an {\em explicit} mechanism has been proposed by Mink and Seifert~\cite{mink-2006-reputation}. There the authors not only propose, but empirically justify a growth rate of
\begin{equation}
\label{eq:mink-seifert}
 \begin{aligned}
 h(R) & = A +  C   \left( 1- \frac{1}{\ln{(e+R)}} \right) \\
         & = A \left[1 +  \frac{C}{A}   \left( 1- \frac{1}{\ln{(e+R)}} \right) \right] \,,
 \end{aligned}
\end{equation}
where $A$ relates to the inherent value of the object for sale and $C$ is a parameter to be fitted. This is accomplished by obtaining data using an auction robot and then computing a single regression, which gives $C=2.50$ in Equation~(\ref{eq:mink-seifert}). To the best of our knowledge, the model in~\cite{mink-2006-reputation} is among the first to give an explicit relationship between reputation and price. We also note that the parameter $\frac{C}{A}$ represents the maximal reputational effect of a seller as a fraction of inherent value $A$, and is what we choose in the numerical examples below. The model also suggests a multiple regression formula where other factors, such as shipping costs and whether a ``buy-it-now" price is offered, are considered as well, in which case $C=1.93$. The authors in~\cite{mink-2006-reputation} also comment that highly experienced sellers have higher feedback scores and design the auction more favorably, which reflects to their higher revenue.
In fact, they show that the coefficient attributed to shipping costs is larger than one, implying that customers put a high value on shipping when deciding on their bids, and that savvy agents take this into consideration. Notice that since $h$ is bounded on $\mathcal{O}$, our admissible class reduces to

\begin{equation}
\mathcal{A}_{\nu} :=  \left \{ \mu \in \mathcal{F}_{s} \mid \mu \text{ is } \mathcal{F}_{s} \text{-p.m.},
\mu_{.} \in U \text{ on } [0,\infty)   \right \}.
\end{equation}

Finally, the work in~\cite{mink-2006-reputation} posits that the horizon does not affect the revenue stream as much as the shipping cost and reputation factors, and so we consider an infinite horizon model here.

\subsection{Hamilton-Jacobi-Bellman Formulation}
\label{sec:HJB}

With the stochastic dynamics for reputation proposed in \cite{sethi-1983-deterministic}, a growth rate model for reputational effect on sales~\cite{mink-2006-reputation}, and an infinite horizon, we expect that switching would depend only on the current reputational state. We now seek a twice-continuously differentiable, polynomially growing function $V \in C^{2}[0,\infty) \cap C_{p}[0,\infty)$ which is a candidate solution of the optimal control problem
\begin{equation}
 \begin{aligned}
 \label{eqn:Opt_Stop}
\bar{V}(R) & = \sup_{\nu} \sup_{\mu \in \mathcal{A}_{\nu}} \mathbb{E} \left[ \int_{0}^{\tau_{0}}  e^{-\rho s}(1-\mu_{s}) h(R^{(\mu)}_{s}) \trds \mid R_{0}=R \right] \\
                  & = \sup_{\nu} \sup_{\mu \in \mathcal{A}_{\nu}} \mathbb{E} \left[ \int_{0}^{\tau_{\mathcal{O}}}  e^{-\rho s}(1-\mu_{s}) h(R^{(\mu)}_{s}) \trds + 1_{\left \{ \tau_{\mathcal{O}} < \infty \right \}} e^{-\rho \tau_\mathcal{O}}   \frac{1+\epsilon}{\rho} h(\infty) 1_{\left \{R^{(\mu)}_{\tau_{\mathcal{O}}} = \infty \right \}}  \mid R_{0}=R \right] \\
h(R) & = A + C   \left( 1- \frac{1}{\ln{(e+R)}} \right)  \\
\tau_{0} & := \inf \left \{t \geq 0 \mid R_{t} = 0 \right \} \,.
\end{aligned}
\end{equation}
One approach to finding $\bar{V}$ would be to solve Equation \eqref{eqn:Opt_Stop} directly. For example, by the definition of $\bar{V}$, it follows that
\begin{equation}
\begin{aligned}
\bar{V}(0) & = 0 \\
\bar{V}(\infty) & = \sup_{\nu} \sup_{\mu \in \mathcal{A}_{\nu}} \mathbb{E} \left[ \int_{0}^{\infty}  e^{-\rho s}(1-\mu_{s}) h(R^{(\mu)}_{s}) \trds \mid R_{0}=\infty \right] \\
& =  \frac{1+\epsilon}{\rho} h(\infty)  \,.
\end{aligned}
\end{equation}

However, we shall instead apply the principle of optimality, and in doing so arrive at the following nonlinear Hamilton-Jacobi-Bellman (HJB) boundary value problem (suppressing the explicit dependence of $R$ on $\mu$)
\begin{equation}
\begin{aligned}
0 & = \max_{-\epsilon \leq \mu \leq \epsilon } \left[ (\mu - \kappa R) \frac{\partial V}{\partial R} + \frac{\sigma^{2}}{2}  \frac{\partial^{2} V}{\partial R^{2}}  + (1-\mu) h(R) - \rho V \right] \\
V(0) & = 0 \\
V(\infty) & = \frac{1+\epsilon}{\rho} h(\infty)  \\
h(R) & = A + C   \left( 1- \frac{1}{\ln{(e+R)}} \right)  =  A + C  \frac{\ln\left(1+\frac{R}{e}\right)}{1+\ln\left(1+\frac{R}{e}\right)}   \,,
\end{aligned}
\end{equation}
which simplifies to
\begin{equation}
\begin{aligned}
\label{eqn:HJB_R}
\rho V  & = h(R) - \kappa R \frac{\partial V}{\partial R} + \frac{\sigma^{2}}{2}  \frac{\partial^{2} V}{\partial R^{2}} +  \epsilon  \Big | \frac{\partial V}{\partial R}  - h(R) \Big | \\
V(0)  & = 0 \\
V(\infty)  & = \frac{1+\epsilon}{\rho} (A+C)  \\
h(R) & = A + C  \frac{\ln\left(1+\frac{R}{e}\right)}{1+\ln\left(1+\frac{R}{e}\right)} \,.
\end{aligned}
\end{equation}
We therefore instead solve the HJB problem, which is justified by the following standard theorem (Corollary IV.5.1 in~\cite{book-fleming-soner-2006}):

\begin{theorem} Let $V \in C^{2}(\mathcal{O}) \cap C_{b}(\bar{\mathcal{O}})$ be a twice-continuously differentiable and bounded solution to an associated HJB equation of a control problem 

\begin{equation}
V_{\rm PM}(x) := \inf_{\nu} \inf_{A_{\nu} } \mathbb{E} \Big[ \int_{0}^{\tau_{\mathcal{O}}} e^{-\rho s} G(R^{\mu}_{s},\mu_{s}) \trds  \Big]
\end{equation}

for a process $R$ whose SDE has drift and volatility that are Lipschitz in $\mu$ and $R$, and a $G(R,\mu)$ that is continuous and polynomially growing for all  $R \in \mathcal{O}$ and continuous for all $\mu \in U$. Moreover, assume either that $\beta > 0$ or that $\tau_{\mathcal{O}} < \infty$ with probability $1$ for every admissible progressively measurable control process $\mu$. Then $V(x)= V_{{\rm PM}}(x)$  for an optimal control $u^{*}(s) \in \argmin [ L_{R}[V] + G(R,u) ]$, where $L$ is the generator of the process $R$.
\end{theorem}

We now have our verification theorem:

\begin{theorem} \label{Verification}
Let $V \in C^{2}(\mathcal{O}) \cap C_{b}(\bar{\mathcal{O}})$ be a twice-continuously differentiable and bounded solution to the boundary value problem Equation \eqref{eqn:HJB_R}. Then $V$ is an optimal solution to the optimal control problem Equation \eqref{eqn:Opt_Stop}.
\end{theorem}
\begin{proof}
By Corollary IV.5.1 in~\cite{book-fleming-soner-2006}, since
  \begin{itemize}
   \item[$(i)$] the drift and volatility are (uniformly) Lipschitz in $R$ and $\mu$,
   \item[$(ii)$] $\rho > 0$,  and
   \item[$(iii)$] we impose an additional boundary condition on the solution as $R \rightarrow \infty$ as our indefinite horizon is now $[0,\tau_{\mathcal{O}}) = [0,\tau_{0}) = [0,\infty)$,
\end{itemize}
it follows directly that our bounded classical solution $V = \bar{V}$ and the control
\begin{equation}
\mu^{*}(R) \in \argmax_{\mu \in [-\epsilon,\epsilon]} \left[ (\mu - \kappa R) \frac{\partial V}{\partial R} + \frac{\sigma^{2}}{2}  \frac{\partial^{2} V}{\partial R^{2}}  + (1-\mu) h(R) - \rho V \right]
\end{equation}
leads to an optimal, stationary Markov control $\mu^{*}(R^{(\mu^{*})}_{s})$.
\end{proof}

\begin{remark} We note here that it is sufficient to require polynomial growth on $V$, not bounded growth. However, as $h$ is bounded above by $A+C$ on $\mathcal{O}$, we can restrict our attention to bounded $V$. Furthermore, an upper bound (and in fact limit as $R \rightarrow \infty$) for $V$ is the value obtained by a seller earning maximum value $A+C$ on her items on the time interval $[0,\infty)$ while always processing orders.
\end{remark}

\section{Market Share Based Pricing Mechanism}
\label{sec:mod.market.share}

We observe that under its current definition, the reputation $R \in \mathcal{O} $. In this section, we employ the~mapping
\begin{equation}
    \label{eqn:map}
	Y := f(R) =  \frac{R}{1+R} \,,
\end{equation}
which results in a normalized \textit{market share} reputation $Y \in \mathcal{Q} := (0,1)$. This is consistent with prior assumptions \cite{mink-2006-reputation} that the value function $h(R)$ be concave, and bounded. We therefore also simplify \eqref{eq:mink-seifert} accordingly, replacing logarithms with a rational function
\begin{equation}
 \begin{aligned}
	\bar{h}(R) & = A + C \frac{R}{1+R} \\
                             & = A \left[ 1 + \frac{C}{A} \frac{R}{1+R} \right] \,.
 \end{aligned}
\end{equation}
This modification, when combined with the mapping from a reputational score to a market share score, $R \rightarrow Y$, produces a more intuitive growth rate, $\tilde{h}(Y)$ which grows linearly with respect to market share reputation.
Recall that $\infty$ is an absorbing state for $R$; consequently, we have $dY_{t} \equiv 0$ if $Y_{t} \geq 1$. For $Y_{t} <1$, it follows from Ito calculus that $dY_{t} = f^{\prime}(R) dR_{t} + \frac{1}{2}f^{\prime \prime}(R) dR_{t} dR_{t}$ and so
\begin{equation}
 \begin{aligned}
dY_{t}  & =  (1-Y)^2 \left[ \Big( \mu - \kappa \frac{Y}{1-Y} \Big) \trdt+ \sigma \trdBt \right] + \frac{1}{2} (-2 (1-Y)^{3}) \sigma^{2} \trdt  \\
             & = \left[ \mu  (1-Y)^2  - \kappa Y(1-Y)   - \sigma^{2}  (1-Y)^{3}  \right]  \trdt  + \sigma  (1-Y)^2   \trdBt  \\
      \tilde{h}(Y)        & = A + CY  \\
\tilde{V}(y) & = \sup_{\nu} \sup_{\mu \in \mathcal{A}_{\nu}} \mathbb{E} \left[ \int_{0}^{\tau_{0}}  e^{-\rho s}(1-\mu_{s}) \tilde{h}(Y_{s}) \trds \mid Y_{0}=y \right] \,.
\end{aligned}
\end{equation}

Note that a function that is bounded for market share $Y \in (0,1)$ is also correspondingly bounded for all reputation  $R \in (0,\infty)$. We therefore study the market share model below, confident that our results will hold for the full reputation model by virtue of the inverse map of Equation \eqref{eqn:map}.
\subsection{Rescaled HJB Model}

After market share rescaling, the HJB now takes a form which is in fact degenerate at both endpoints  $y=0,1$:
\begin{equation}
\begin{aligned}
\label{eqn:HJB}
	&- \frac{\sigma^{2}}{2} (1-y)^{4} V^{\prime \prime}(y) +   \rho V  =  \tilde{h}(y)  - \left[ \kappa y(1-y) +  \sigma^{2}  (1-y)^{3} \right] V^{\prime}(y)   +  \epsilon  \Big | (1-y)^{2} V^{\prime}(y)  - \tilde{h}(y) \Big | \\
	&V(0)  = 0, \qquad V(1)  = \frac{1+\epsilon}{\rho}\tilde{h}(1) \,.
\end{aligned}
\end{equation}
In theory, a closed form solution for this boundary value problem can be formally constructed piecewise, to the left and right of the special point $y^{*}$, for which $\tilde{h}(y^*) = (1-y^*)^{2} V^{\prime}(y^*)$. We shall refer to $y^*$ as the switching point below, since it is the reputational value at which the seller switches from advertising to processing. In practice, we shall resort to numerical computation to approximate this solution, and in doing so estimate the value $y^*$.

\begin{theorem}
If $V^{**} \in C^2(\mathcal{Q})$ is a monotonically increasing solution to the HJB problem Equation \eqref{eqn:HJB}, then  $V^{**} = \tilde{V}$.
\end{theorem}

\begin{proof}
 Enforcing the boundary condition on $V^{**}(1)$ enforces bounded growth on our monotonically increasing $V^{**}$. By the inverse map of Equation \eqref{eqn:map}, the result follows as a corollary to Theorem~\ref{Verification}. Note that $\mathcal{Q} = (0,1)$ is open, but $y=1$ is not reached in finite time from anywhere in $\mathcal{Q}$ and the diffusion $Y$ is in fact absorbed for all time in state $Y=1$. Hence, the exact value of $V^{**}(1)$ imposed.
\end{proof}

\subsection{Reduced Model}
\label{sec:red.mod}

With the previous mapping, note that the market (reputation) share $Y$ is absorbed at $Y=1$, which corresponds to $R \rightarrow \infty$. Moreover, as in the Stochastic Nerlove evolution, we have that the probability that $Y_{t} < 0$ for some positive $t$ is greater than 0, and the drift term is a third-order polynomial in $Y$. Based on these observations, we propose a reduced model for {\it Stochastic Reputation Share}, expressed by the following stochastic differential equation and associated stochastic optimal control problem
\begin{equation}
 \begin{aligned}
  \check{V}(y) & = \sup_{\nu} \sup_{\mu \in \mathcal{A}_{\nu}} \mathbb{E} \left[ \int_{0}^{\tau_{0}}  e^{-\rho s}(1-\mu_{s}) (A + CY_{s}) \trds \mid Y_{0}=y \right] \\
 \frac{dY_{t}}{1-Y_{t}} & =  \mu \trdt +  \sigma \trdBt  \,.
 \end{aligned}
\end{equation}
This leads to a corresponding nonlinear HJB boundary value problem,
\begin{equation}
\begin{aligned}
0 & = \max_{-\epsilon \leq \mu \leq \epsilon } \left[ \mu(1-y) V^{\prime}(y) + \frac{\sigma^{2}}{2} (1-y)^{2} V^{\prime \prime}(y)  + (1-\mu) (A+Cy) - \rho V \right] \\
V(0) & = 0 \\
V(1) & = \sup_{\nu} \sup_{\mu \in \mathcal{A}_{\nu}} \mathbb{E} \left[ \int_{0}^{\infty}  e^{-\rho s}(1-\mu_{s}) \tilde{h}(1) \trds  \right] =
\rho^{-1}(1+\epsilon) \tilde{h}(1)\,,
\end{aligned}
\end{equation}
which can be recast as
\begin{equation}
\begin{aligned}
\label{eqn:HJB_Ry}
-\frac{\sigma^{2}}{2} (1-y)^{2} V^{\prime \prime}(y) + \rho V &= \tilde{h}(y)  +  \epsilon  \Big | (1-y) V^{\prime}(y)  - \tilde{h}(y) \Big | \\
V(0) &= 0 \\
V(1) &= \frac{1+\epsilon}{\rho}\tilde{h}(1) \,.
\end{aligned}
\end{equation}

One of the most attractive features of the reduced model \eqref{eqn:HJB_Ry} is that it has a closed form, piecewise defined analytic solution, which we derive in the Appendix. Here we state the following theorem.
\begin{theorem}
There exists a solution $V^{*} \in C^2(\mathcal{Q}) \cap C_{b}(\bar{\mathcal{Q}})$ to the reduced model \eqref{eqn:HJB_Ry}. Furthermore, $V^{*} = \check{V}$, and is given by the piecewise solution
\begin{align}
	V^{*}(y) =
	\begin{cases}
	V_\ell(y),	&	0\leq y \leq y^*\\
	V_r(y),	&	y^*\leq y \leq 1 \,,
	\end{cases}
\end{align}
where
\begin{align}
	V_\ell &= c_1(1-y)^{\gamma^\ell_-}+c_2(1-y)^{\gamma^\ell_+} + \alpha_\ell + \beta_\ell y \\
	V_r &= c_3(1-y)^{\gamma^r_+} + \alpha_r + \beta_r y \,,
\end{align}
with constants defined in the Appendix.
\end{theorem}
\begin{proof}
 The construction of this piecewise solution is found in the Appendix. The solution constructed is monotonically increasing in $y$. Enforcing the boundary condition on $V^{*}(1)$ uniformly bounds $V^{*}$ on $\mathcal{Q}$. The proof of equality $V^{*} = \check{V}$ then follows from Corollary IV.5.1 in~\cite{book-fleming-soner-2006}, as in our Theorem $2$ above.
\end{proof}
\begin{remark}
Since the reduced model has a known closed form solution, we can measure the error made in constructing numerical solutions in Section \ref{sec:num.results}, which in turn acts as a benchmark for the code we use to study the full model. It is uncommon to find a closed form solution for most optimal control problems.
\end{remark}

\section{Numerical Results}
\label{sec:num.results}
In this section, we construct numerical solutions of the market share scaled model~Equation\eqref{eqn:HJB}, as well as our proposed reduced model~ Equation \eqref{eqn:HJB_R}. In Sections ~\ref{sec:sol.ode} and~\ref{sec:con.solution}, the piecewise analytic solution and a nonlinear equation for the reduced switching point $y^*$ are found. These will be used to validate the numerical results for the reduced model, which act as a benchmark for the full problem. Please note that in what follows, the symbol $y^{*}$ is used to denote the switching point in reputational share for both the reduced and full model and their corresponding ODEs. 

\subsection{Numerical Results for the Reduced Model}
\begin{figure}[!h]
\centering
\subfigure[]{\includegraphics[width = 0.45 \textwidth]{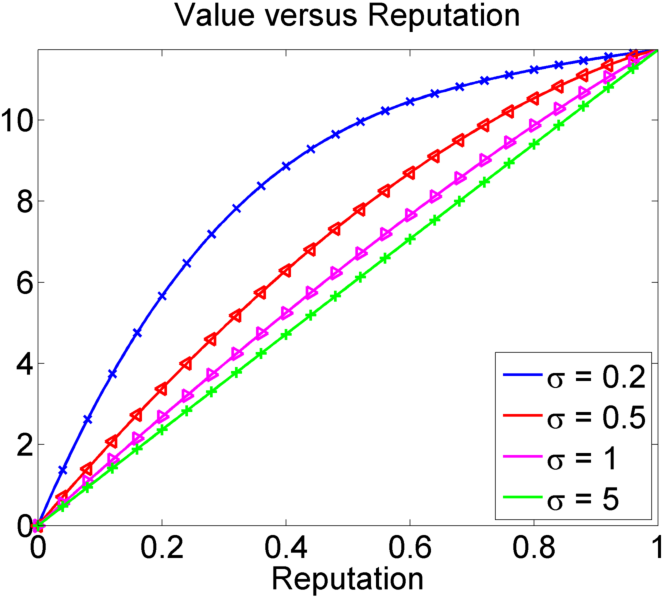}}
\subfigure[]{\includegraphics[width = 0.44 \textwidth]{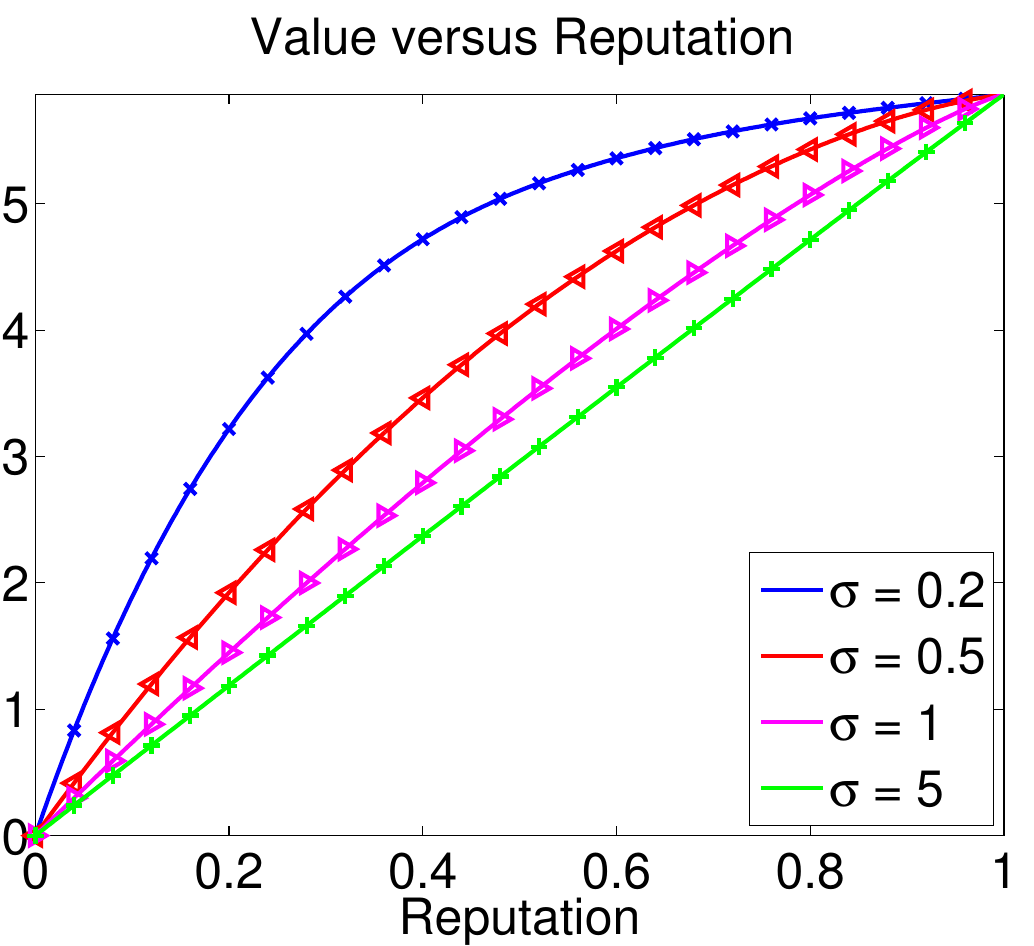}}\\
\subfigure[]{\includegraphics[width = 0.45 \textwidth]{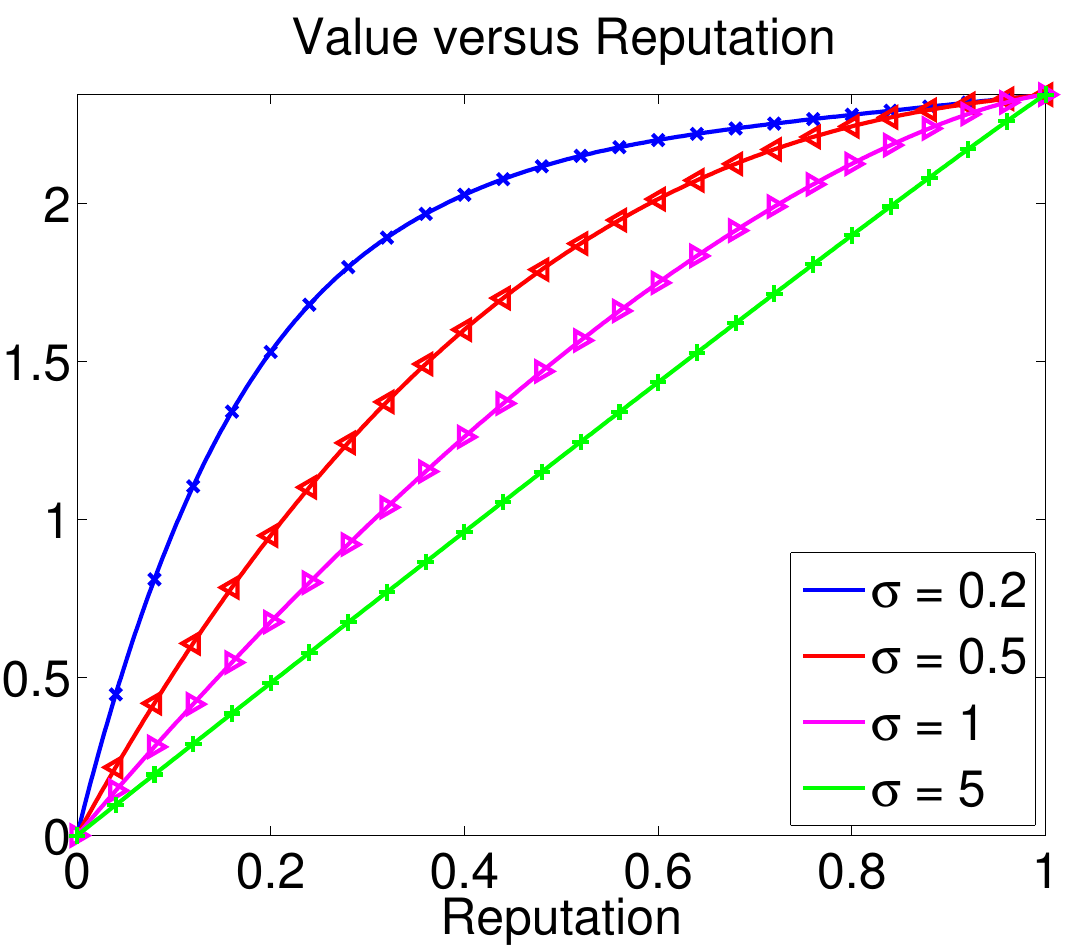}}
\subfigure[]{\includegraphics[width = 0.45 \textwidth]{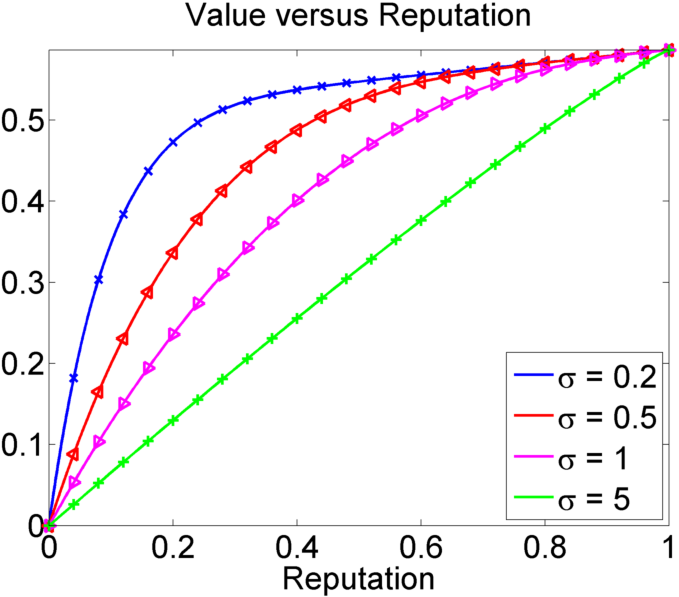}}
\caption{A plot of several numerical solutions $V(y)$ for the reduced boundary value problem Equation \eqref{eqn:HJB_Ry} are shown with $A = 1$, $C = 0.15$, and $\epsilon = 0.02$. ({\bf a}) $\rho = 0.1$; ({\bf b}) $\rho = 0.2$; ({\bf c})~$\rho = 0.5$; ({\bf d}) $\rho = 2.0$.}
\label{fig:Value}
\end{figure}

The boundary value problem Equation \eqref{eqn:HJB_R} is discretized using a standard finite difference scheme. Let $y_k = k/N_y$, and $V_k = V(y_k)$, for $k = 0,1,\ldots, N_y$. Then we solve a linear system of $N_y+1$ equations of the form $M \mathbf{v} = \mathbf{f}$, where
\begin{eqnarray*}
M\mathbf{v} &\approx& \rho V - \frac{\sigma^{2}}{2} (1-y)^{2} V^{\prime \prime} \,, \\ 
\mathbf{f} &\approx& h +  \epsilon  \Big | (1-y) V^{\prime}  - h \Big | \,.
\end{eqnarray*}
Since $\mathbf{f}$ depends on $V$, the solution is implicit, and therefore must be obtained by using a fixed point iteration $M\mathbf{v}^{(k+1)} = \mathbf{f}^{(k)}$. In our numerical experiments, we let $N_y = 1000$, and set the maximum iteration count at $K = 20$. Convergence is observed in all tested cases. In Figure~\ref{fig:Value}, several numerical solutions $V(y)$ are shown, for fixed $A = 1$, $C = 0.15$, $\epsilon = 0.02$, and various values of $\rho=0.1,0.2,0.5,2.0$ and $\sigma=0.2,0.5,1.0,5.0$. The values of $A$ and $C$ are chosen via $\frac{C}{A} = 0.15$ to reflect a maximal $15 \%$ reputational premium for sellers above the inherent value $A$.

\begin{figure}[!h]
\centering
\subfigure[]{\includegraphics[width = 0.45 \textwidth]{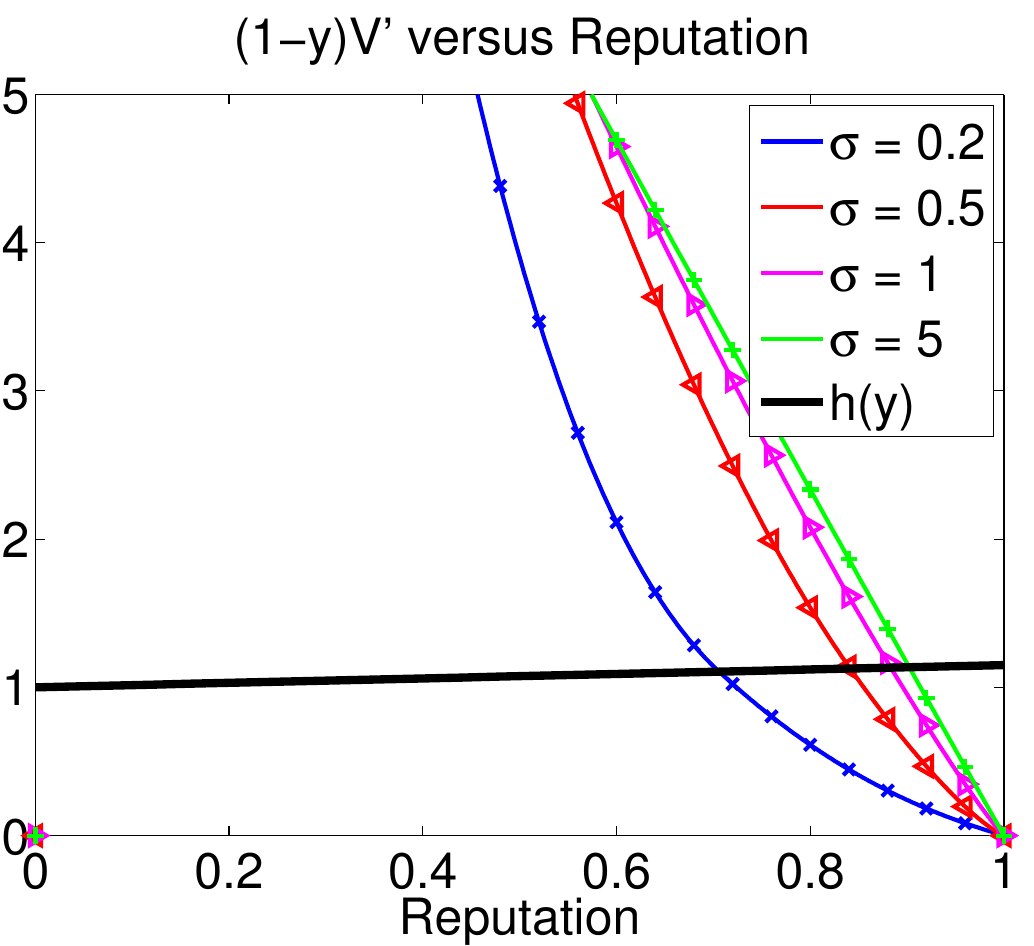}}
\subfigure[]{\includegraphics[width = 0.45 \textwidth]{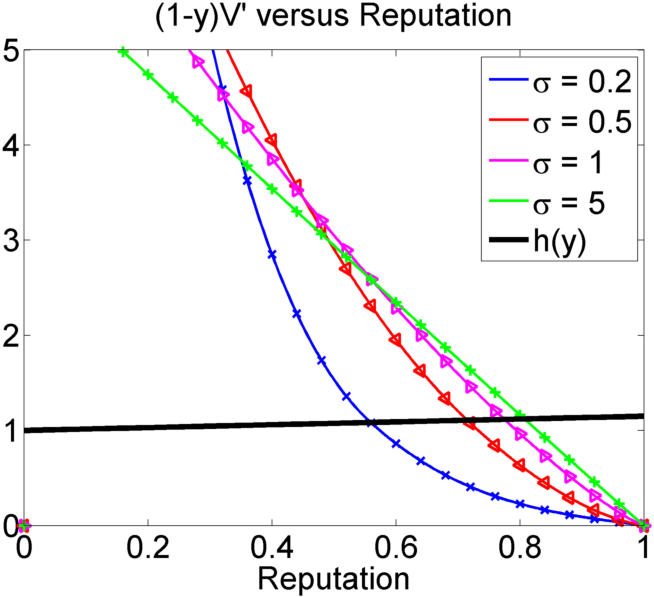}}\\
\subfigure[]{\includegraphics[width = 0.45 \textwidth]{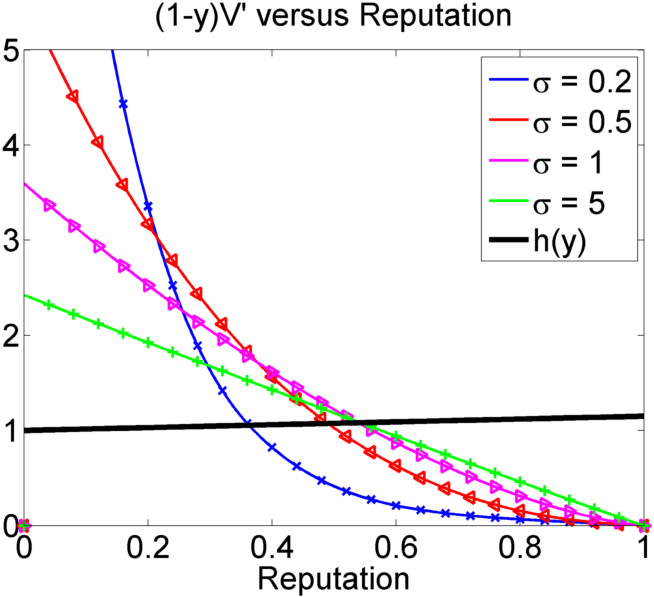}}
\subfigure[]{\includegraphics[width = 0.45 \textwidth]{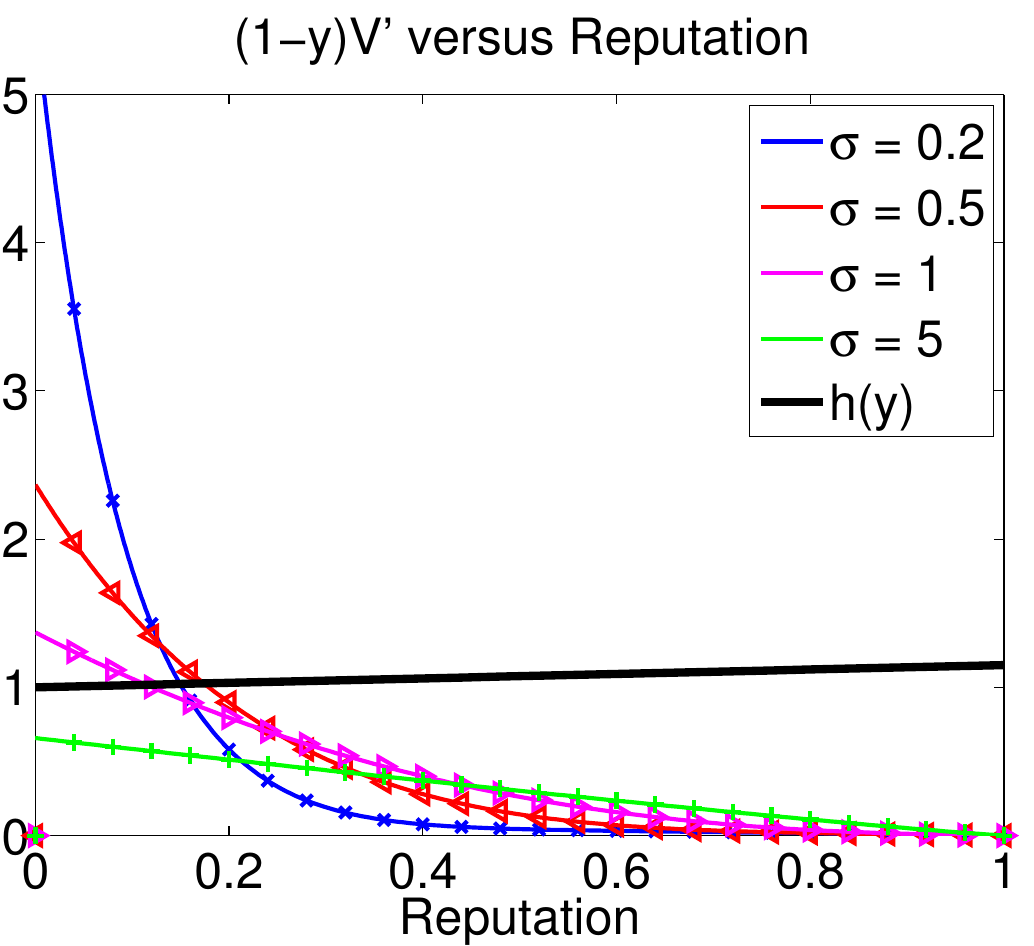}}
\caption{A plot of the quantity $(1-y)V'(y)$, where $V(y)$ is a numerical solution of the the reduced boundary value problem Equation \eqref{eqn:HJB_Ry}. The value $y^*$ is defined as the intersection of these curves with $h(y)$ (the solid line). ({\bf a}) $\rho = 0.1$; ({\bf b}) $\rho = 0.2$; ({\bf c}) $\rho = 0.5$; ({\bf d}) $\rho = 2.0$.}
\label{fig:Switch}
\end{figure}

It follows that both $\rho$ and $\sigma$ have a dramatic effect not only the shape of the solution, but also, as is demonstrated in Figure~\ref{fig:Switch}, on the switching point $y^*$. These numerical results suggest that both the discount rate $\rho$ and the volatility $\sigma$ of a seller's reputation can ``dramatically'' affect her behavior, particularly the critical reputation at which she switches from the processing to advertising mode.

\subsection{Numerical Results for the Full Model}

The same numerical discretization is employed to solve the full Nerlove-Arrow model Equation \eqref{eqn:HJB}, where we solve a linear system of the form $M \mathbf{v} = \mathbf{f}$, with
\begin{eqnarray*}
M\mathbf{v} &\approx& \rho V + \left[ \kappa y(1-y) +  \sigma^{2}  (1-y)^{3} \right] V^{\prime}(y) - \frac{\sigma^{2}}{2} (1-y)^{4} V^{\prime \prime}(y) \,, \\
\mathbf{f} &\approx& h +  \epsilon  \Big | (1-y)^2 V^{\prime}  - h \Big | \,.
\end{eqnarray*}

\begin{figure}[!h]
\centering
\subfigure[]{\includegraphics[width = 0.45 \textwidth]{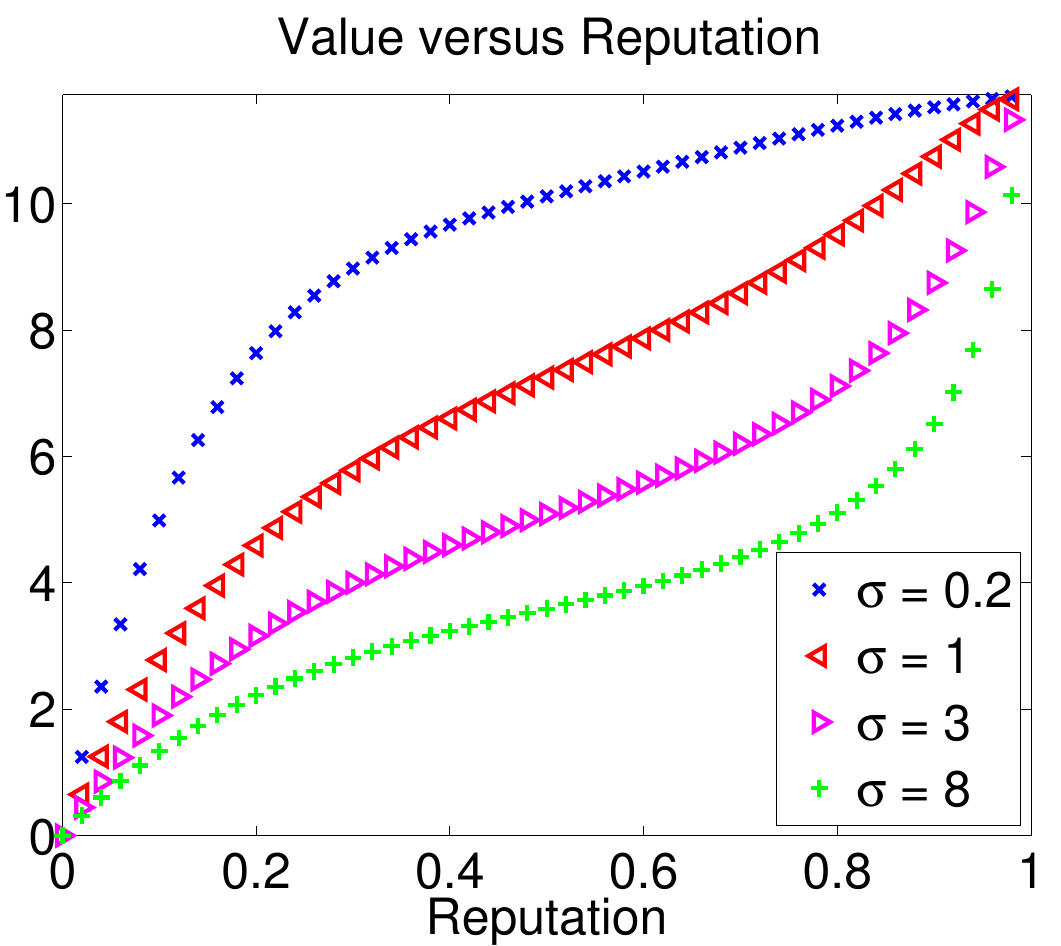}}
\subfigure[]{\includegraphics[width = 0.44 \textwidth]{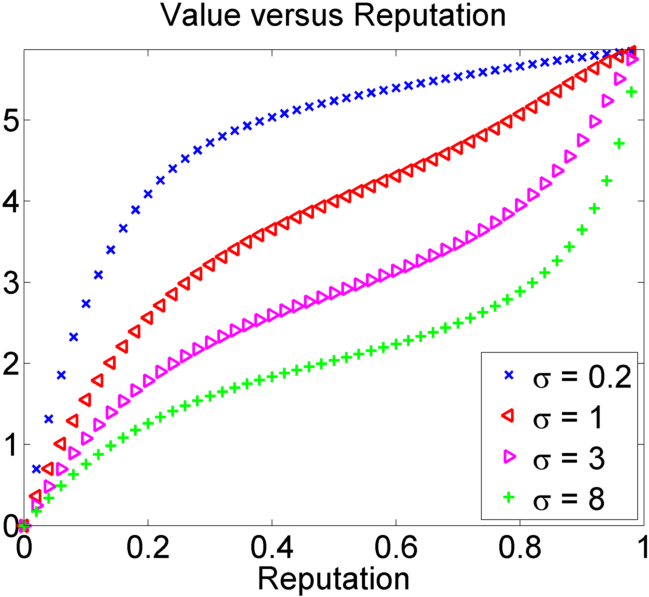}}\\
\subfigure[]{\includegraphics[width = 0.45 \textwidth]{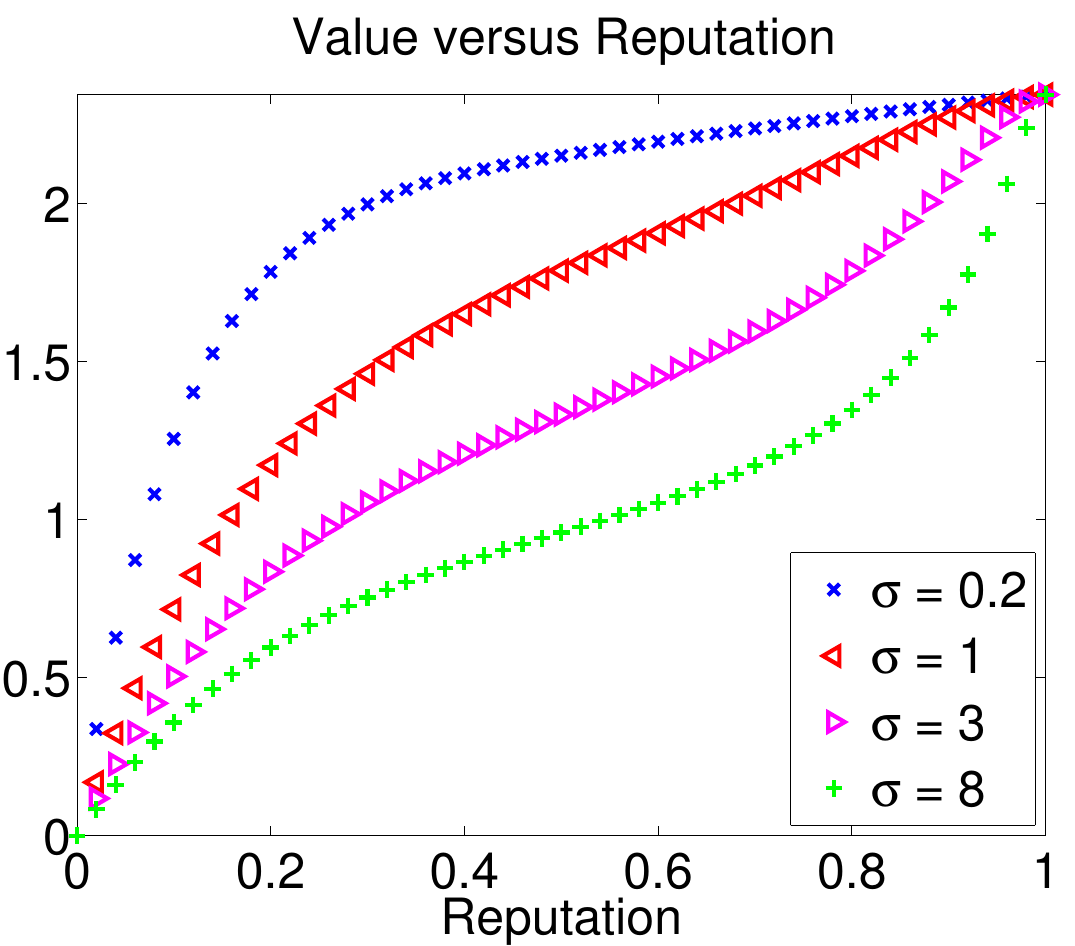}}
\subfigure[]{\includegraphics[width = 0.45 \textwidth]{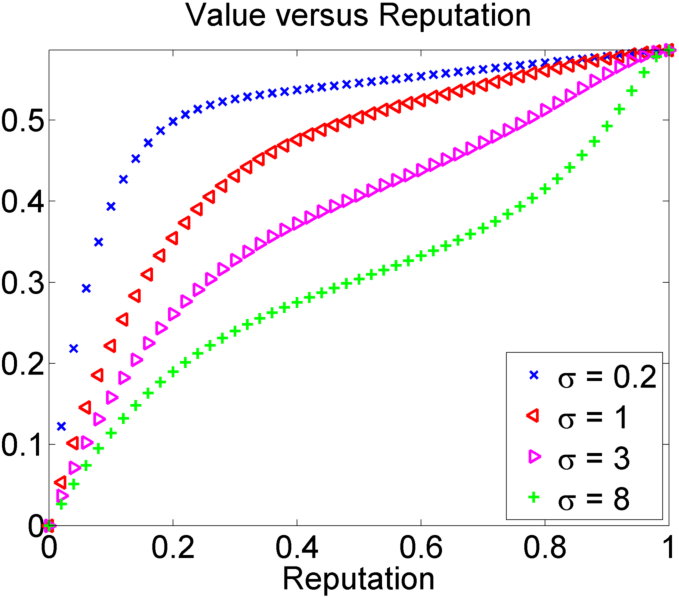}}
\caption{A plot of several numerical solutions $V(y)$ for the full boundary value problem Equation \eqref{eqn:HJB} are shown with $A = 1$, $C = 0.15$, $\epsilon = 0.02$, and $\kappa = 0$. ({\bf a}) $\rho = 0.1$; ({\bf b})~$\rho = 0.2$; ({\bf c}) $\rho = 0.5$; ({\bf d}) $\rho = 2.0$.}
\label{fig:Full_Value}
\end{figure}

\begin{figure}[!h]
\centering
\subfigure[]{\includegraphics[width = 0.45 \textwidth]{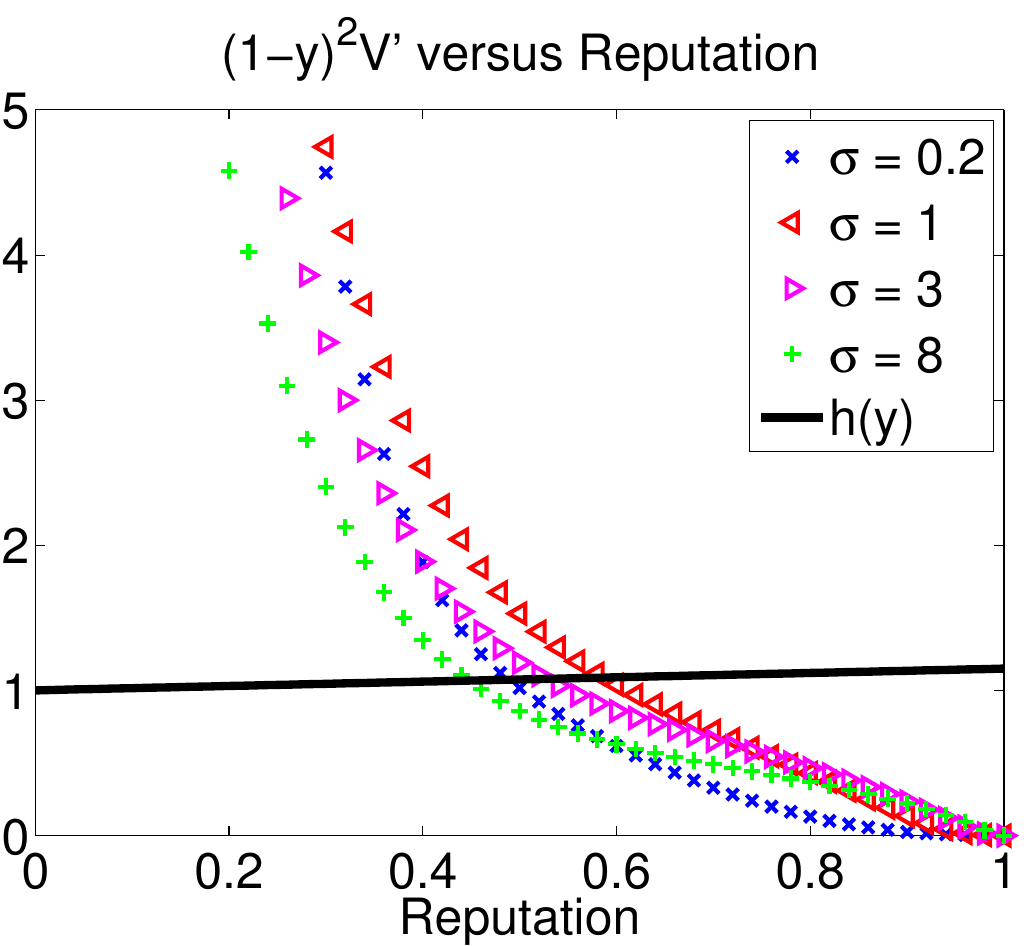}}
\subfigure[]{\includegraphics[width = 0.45 \textwidth]{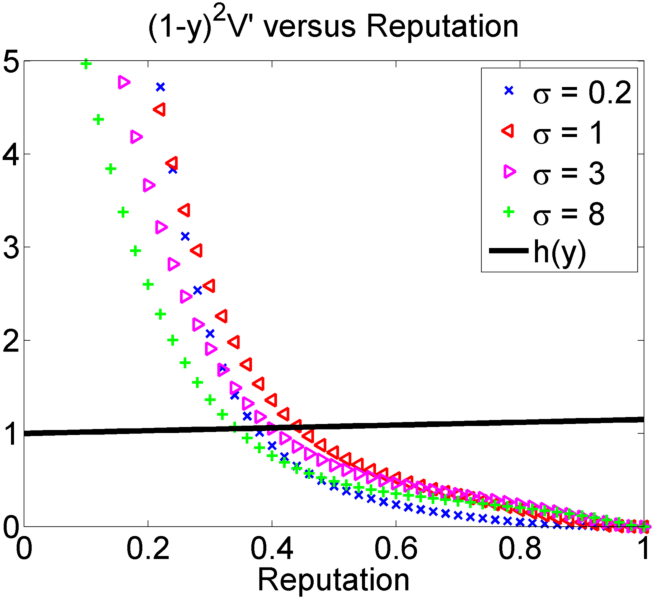}}\\
\subfigure[]{\includegraphics[width = 0.45 \textwidth]{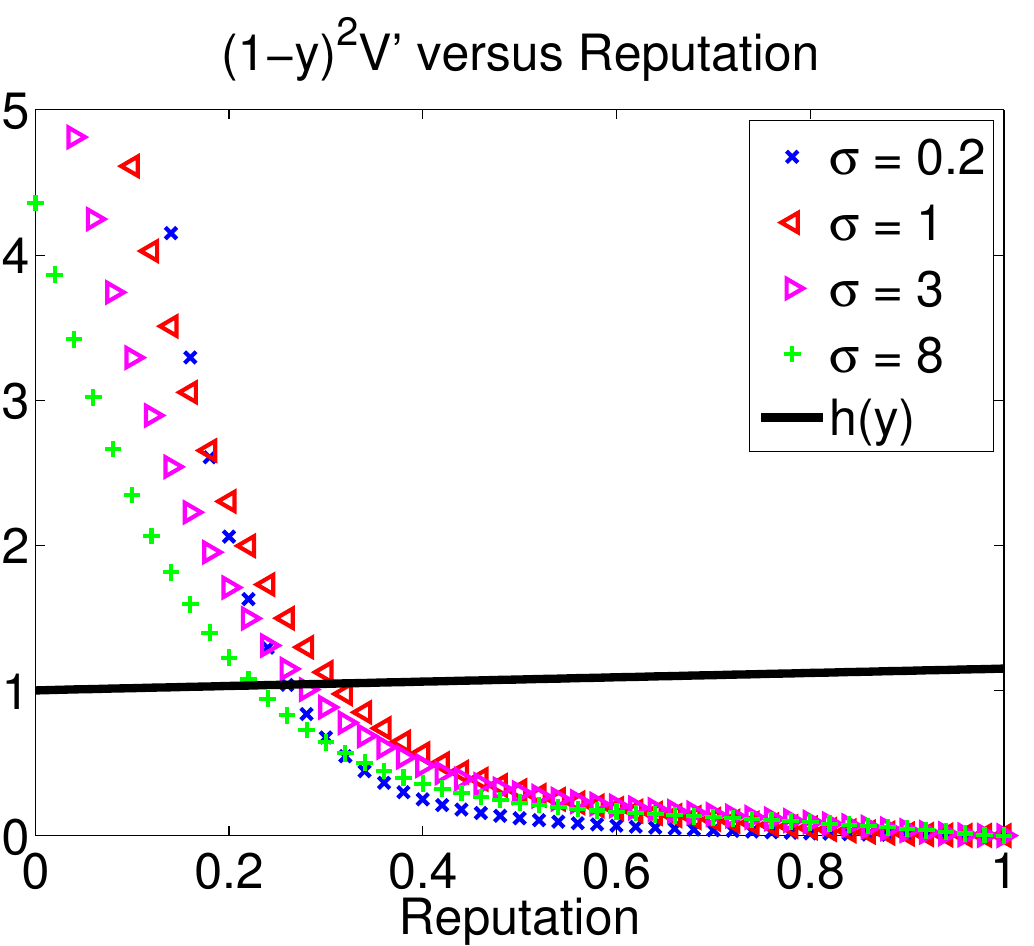}}
\subfigure[]{\includegraphics[width = 0.45 \textwidth]{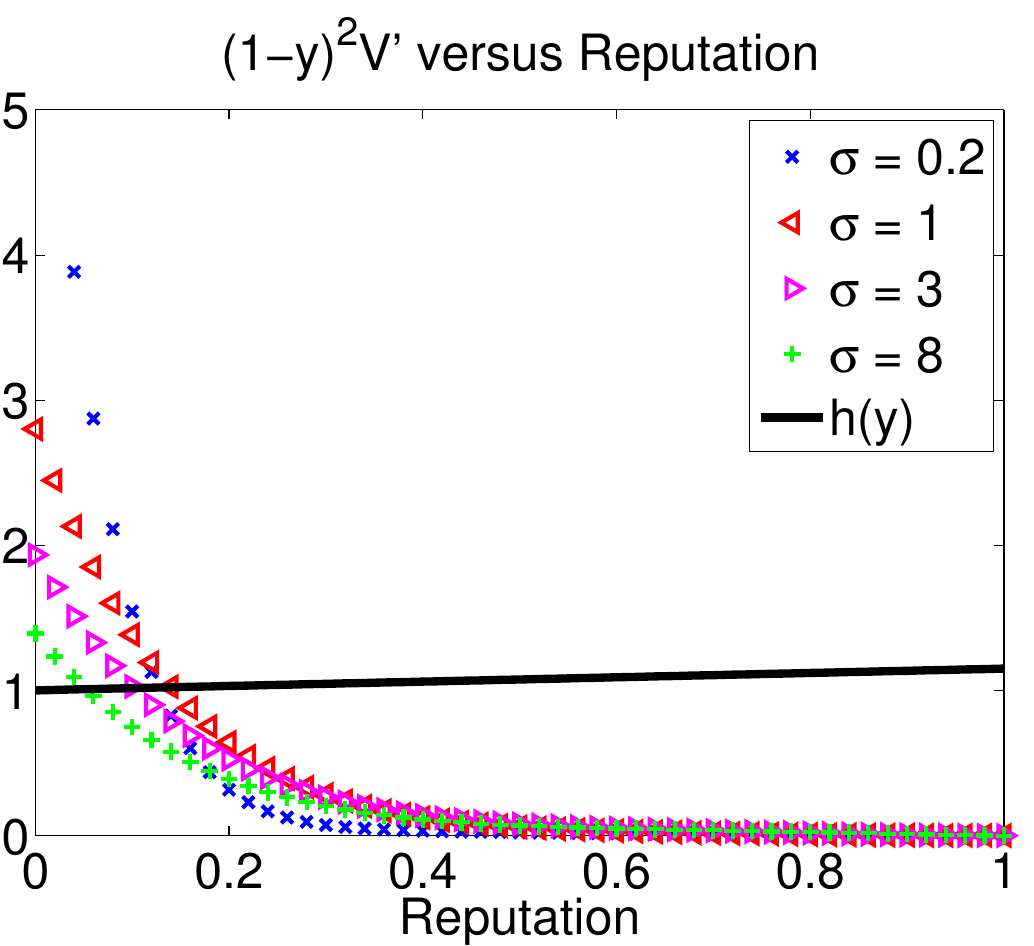}}
\caption{A plot of the quantity $(1-y)^2V'$, where $V$ is a numerical solution of the the full boundary value problem Equation \eqref{eqn:HJB}. The value $y^*$ is defined as the intersection of these curves with $h(y)$ (black solid line). ({\bf a}) $\rho = 0.1$; ({\bf b}) $\rho = 0.2$; ({\bf c}) $\rho = 0.5$; ({\bf d}) $\rho = 2.0$.}
\label{fig:Full_Switch}
\end{figure}

We first set $\kappa = 0$, and hold all remaining parameters fixed, and plot the results in Figures~\ref{fig:Full_Value} and~\ref{fig:Full_Switch}. The full solutions have more curvature than those of the reduced model, but nonetheless remain monotone, and have a single unique switching point $y^*$. The additional curvature is due to the appearance of $V'$ terms in the differential operator, which now depend on $\sigma$, as well as $\kappa$, which we recall incorporates mean reversion. That is, independent of the sellers strategy, reputation will tend to decrease to a smaller amount of the market share, with constant rate $\kappa$.

\begin{figure}[!h]
\centering
\subfigure[]{\includegraphics[width = 0.45 \textwidth]{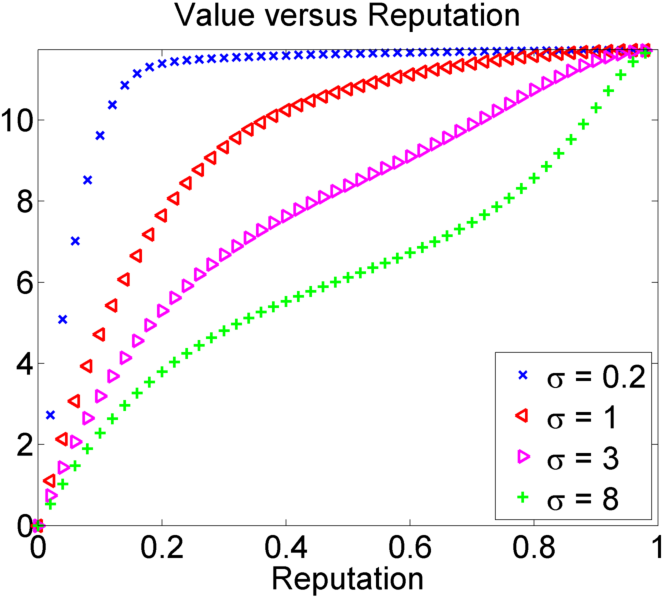}}
\subfigure[]{\includegraphics[width = 0.44 \textwidth]{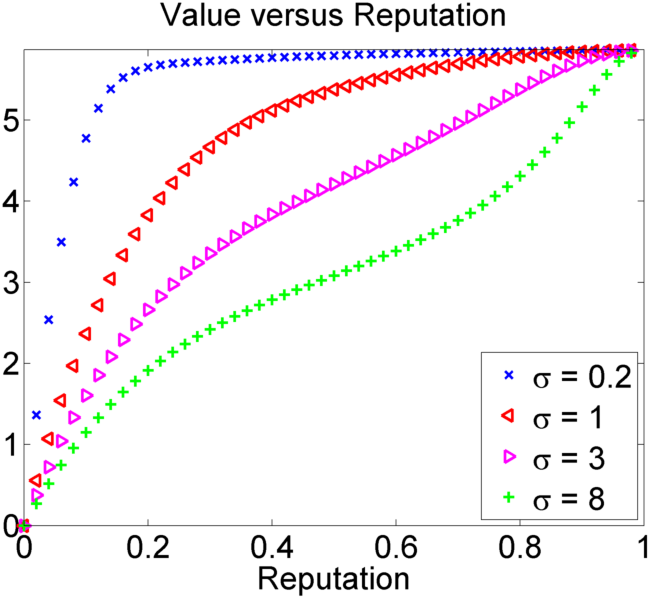}}\\
\subfigure[]{\includegraphics[width = 0.45 \textwidth]{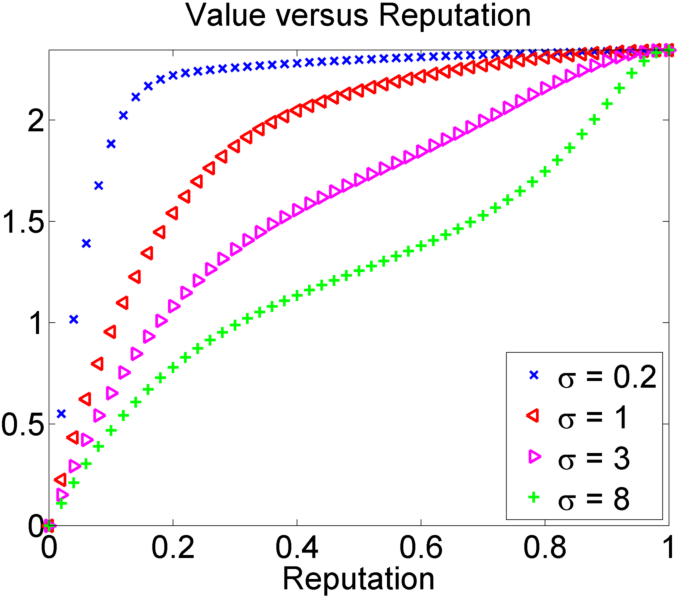}}
\subfigure[]{\includegraphics[width = 0.45 \textwidth]{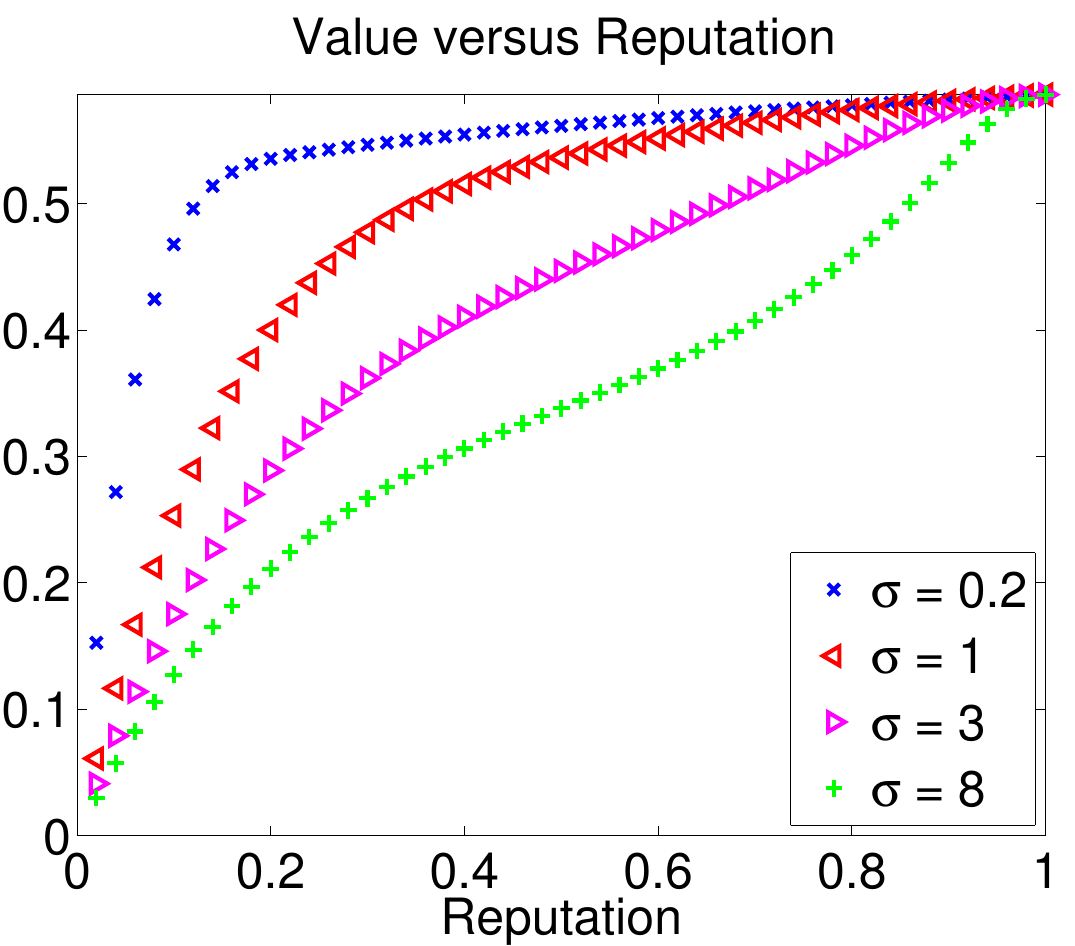}}
\caption{The same as Figure~\ref{fig:Full_Value}, but with $\kappa = 1$. ({\bf a}) $\rho = 0.1$; ({\bf b}) $\rho = 0.2$; ({\bf c}) $\rho = 0.5$; ({\bf d})~$\rho = 2.0$.}
\label{fig:Full_Value_2}
\end{figure}

\begin{figure}[!h]
\centering
\subfigure[]{\includegraphics[width = 0.45 \textwidth]{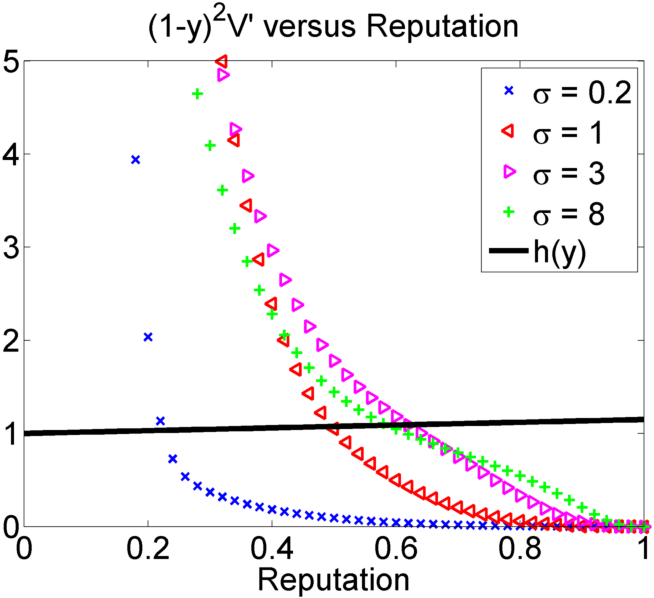}}
\subfigure[]{\includegraphics[width = 0.45 \textwidth]{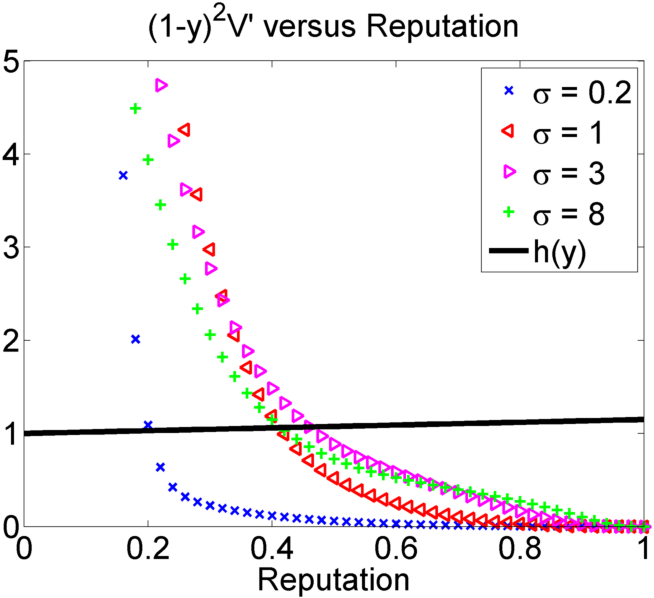}}\\
\subfigure[]{\includegraphics[width = 0.45 \textwidth]{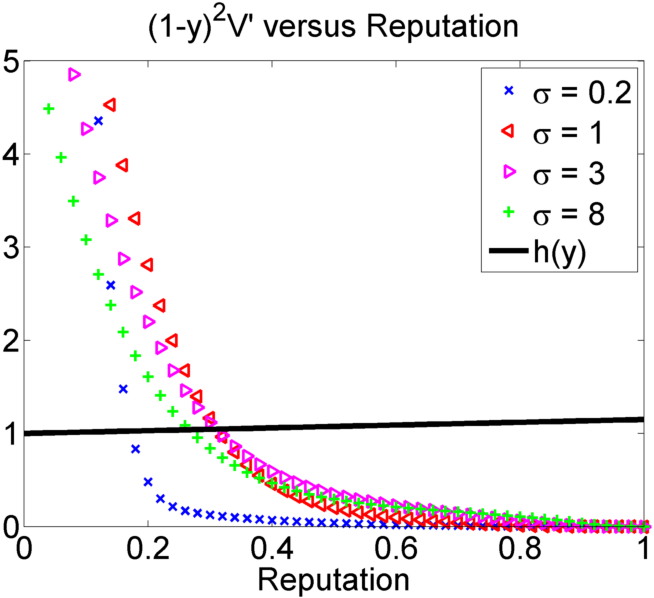}}
\subfigure[]{\includegraphics[width = 0.45 \textwidth]{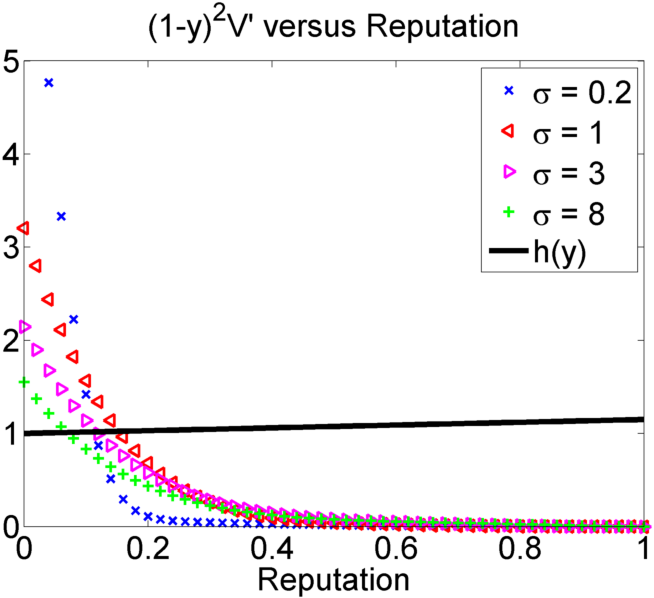}}
\caption{The same as Figure~\ref{fig:Full_Switch}, but with $\kappa = 1$. ({\bf a}) $\rho = 0.1$; ({\bf b}) $\rho = 0.2$; ({\bf c}) $\rho = 0.5$; ({\bf d})~$\rho = 2.0$.}
\label{fig:Full_Switch_2}
\end{figure}

In Figures \ref{fig:Full_Value_2} and \ref{fig:Full_Switch_2}, the same solutions are shown with $\kappa = 1$. Here it becomes apparent that both the rate of mean reversion, as well as the volatility will affect the seller's optimal strategy.
\begin{remark}
Our numerical results illustrate that a wide range of seller behaviors can be described by varying the parameters, and therefore that reputational value is a strong indicator of optimal seller~behavior.
\end{remark}

\section{Conclusions and Future Work}
\label{sec:conclusion}

In this work, we have proposed a model which describes the behavior of sellers participating in online auctions (markets). Our model is based on the premise that buyer feedback greatly impacts the sales rate, which motivates our assumption that the revenue per sale (or the price per unit) is solely dependent upon seller reputation. We assumed that a seller has a fixed amount of resources that must be allocated either to (i) advertising to new buyers or (ii) processing orders for current customers. In doing so, we were able to design an optimal selling strategy, wherein the seller switches their behavior when an optimal market share reputation is reached, which depends on statistical modeling parameters.

These modeling parameters were introduced through the wealth-reputation mechanisms in~\cite{nerlove-1962-optimal}, which have been generalized to stochastic settings in~\cite{sethi-1983-deterministic, rao-1986-estimating, raman-2006-boundary}. We have additionally modified the empirical model $h(R)$ found in~\cite{mink-2006-reputation}, relating the price per unit to the reputation. Rather than viewing reputation $R$ as an unbounded quantity, we instead introduced a market scaled reputation $Y = R/(1+R)$, where $Y \in [0,1)$, which reduces the price per unit to a simple linear model $h(Y) = A+CY$.

We then optimized the stochastic model over the control of the excess rate $\mu$ using the Hamilton-Jacobi-Bellman equation, yielding a deterministic equation relating the value per unit to the seller reputation. The resulting boundary value problem is then solved numerically, and we qualitatively see that the value of sellers goods increase monotonically with reputation, but that a unique optimal reputation $y = y^*$ determines when the seller should switch from advertising to processing. The numerical scheme was validated with a reduced model, which has a closed form piecewise analytic solution, and permits direct determination of $y^*$. The numerical results validate our modeling assumptions, and provide a framework for studying seller behavior based on seller reputation.
Although the used techniques are standard, we believe that the optimal strategy presented both analytically and numerically has implications on the way reputational information can be used to predict the behavior of an individual seller in an online market.

This work can naturally be generalized in a variety of ways. For instance, the model could be developed in the case of finite horizon time $T$, since a more realistic assumption is that a seller strategy depends on both time and the current reputation state. Unlike our current infinite horizon model, the resulting HJB equation will lead to a time-dependent partial differential equation, which will be the subject of future work.

\appendix
\section{Piecewise Analytic Solution of Reduced Model}
\label{sec:sol.ode}

We now construct the exact solution to the reduced HJB boundary value problem Equation \eqref{eqn:HJB_R}, which incidentally can be found analytically. The principal difficulty that must be overcome is the appearance of an absolute value, which contains the unknown variable $V'$. We therefore consider a piecewise defined solution $V(y)$, which remains $C^2$ across the switching point $y^*$, defined by the vanishing of the expression in the absolute value
\begin{equation}
\label{eqn:ystar}
(1-y^*)V'(y^*) = \tilde{h}(y^*) \,.
\end{equation}
%
We now separately consider the regions $y<y^*$ and $y>y^*$ in the reduced model Equation \eqref{eqn:HJB_R}, and define the corresponding linear differential operators $L_{\pm}$ as
\begin{equation}
L_{\pm}[V] = \left(-\frac{\sigma^2}{2}(1-y)^2\frac{d^2}{dy^2}\pm\epsilon(1-y)\frac{d}{dy}+\rho\right)V \,.
\end{equation}
The boundary value problem can then be decomposed into two smaller problems for $V_\ell$ and $V_r$.
The boundary conditions provide one condition for each of $V_\ell$ and $V_r$.
The remaining two conditions are provided by enforcing $C^2$ smoothness of the solution at the switching point Equation \eqref{eqn:ystar}.
Hence we now formulate two (well-posed) boundary value problems:
\begin{align}
	\label{eqn:ODEL}
	L_- [V_{\ell}] &= (1-\epsilon)\tilde{h}(y) 	\qquad	 0<y<y^* \\
	\label{eqn:BCL}
	V_{\ell}(0) &= 0,				\qquad	(1-y^*)V_{\ell}^{\prime}(y^*) = \tilde{h}(y^*)
\end{align}
and
\begin{align}
	\label{eqn:ODER}
	L_+ [V_{r}] &= (1+\epsilon)\tilde{h}(y)			\qquad	y^*<y<1\\
	\label{eqn:BCR}
	V_{r}(1) &= \frac{1+\epsilon}{\rho}\tilde{h}(1),	\qquad	(1-y^*)V_r^{\prime}(y^*) = \tilde{h}(y^*) \,.
\end{align}
Once solved, enforcing continuity uniquely defines the value of the switching point $y^*$:
\begin{equation}
V_{\ell}(y^*) = V_r(y^*) \,.
\end{equation}

\subsection{Constructing the Solution}
\label{sec:con.solution}

We shall construct the solutions $V_\ell$ and $V_r$ separately first, although the approach for both solutions will be the same. Following standard methods, we first decompose the full solution into a homogeneous and particular solution, where the homogeneous part satisfies $L_{\pm} [u] = 0$. Since the operators $L_{\pm}$ are equi-dimensional, we seek solutions of the form
\begin{equation}
u = c (1-y)^{\gamma} \,,
\end{equation}
and the application of the differential operator yields
\begin{equation}
L_{\pm} [u] = c(1-y)^\gamma \left(-\frac{\sigma^2}{2}(\gamma^2-\gamma) \mp \epsilon \gamma +\rho\right) = 0 \,.
\end{equation}
We set the parenthetical term of this expression to zero in order to solve for the admissible exponents $\gamma$
\begin{align}
	\gamma^{\ell}_{\pm}	&= \frac{1}{2} + \frac{\epsilon}{\sigma^2}	\pm \sqrt{ \left(\frac{1}{2} + \frac{\epsilon}{\sigma^2}\right)^2+2\frac{\rho}{\sigma^2}} &\text{(from $L_-$)}\\
	\gamma^{r}_{\pm}	&= \frac{1}{2} - \frac{\epsilon}{\sigma^2}	\pm \sqrt{ \left(\frac{1}{2} - \frac{\epsilon}{\sigma^2}\right)^2+2\frac{\rho}{\sigma^2}} &\text{(from $L_+$)} \,.
\end{align}
Because $\gamma^r_- <0$ for $\rho>0$, the homogeneous solution $(1-y)^{\gamma^r_-}$ will be unbounded as $y\to 1$, and so we exclude it from further consideration below. Furthermore, since $\tilde{h}(y) = A + Cy$ is a linear function, the particular solution will also be linear. We therefore have a general solution of the form
\begin{align}
	V_\ell &= c_1(1-y)^{\gamma^\ell_-}+c_2(1-y)^{\gamma^\ell_+} + \alpha_\ell + \beta_\ell y \\
	V_r &= c_3(1-y)^{\gamma^r_+} + \alpha_r + \beta_r y \,.
\end{align}
The particular solution is fixed by enforcing that the differential Equations \eqref{eqn:ODEL} and \eqref{eqn:ODER} are satisfied, which yields
\begin{align}
	\alpha_\ell &= \left(\frac{1-\epsilon}{\rho}\right)(A+C)-\left(\frac{1-\epsilon}{\rho+\epsilon}\right)C, \qquad \beta_\ell = \left(\frac{1-\epsilon}{\rho+\epsilon}\right)C \\
	\alpha_r &= \left(\frac{1+\epsilon}{\rho}\right)(A+C)-\left(\frac{1+\epsilon}{\rho-\epsilon}\right)C, \qquad \beta_r = \left(\frac{1+\epsilon}{\rho-\epsilon}\right)C\,.
\end{align}
The homogeneous coefficients $c_1$ and $c_2$ are then found by applying the boundary conditions \eqref{eqn:BCL}, and the resulting linear system
\[
	\begin{pmatrix}
	1 & 1 \\
	\gamma^\ell_- (1-y^*)^{\gamma^\ell_-} & \gamma^\ell_+ (1-y^*)^{\gamma^\ell_+}
	\end{pmatrix}
	\begin{pmatrix}
	c_1 \\
	c_2
	\end{pmatrix}
	=
	\begin{pmatrix}
	-\alpha_\ell \\
	\beta_\ell(1-y^*) - h(y^*)
	\end{pmatrix}
\]
is solved by
\begin{align}
	\label{eqn:c1_def}
	c_1 &= -\frac{\alpha_\ell\gamma^\ell_+ (1-y^*)^{\gamma^\ell_+} +\beta_\ell(1-y^*) - \tilde{h}(y^*) }{\gamma^\ell_+ (1-y^*)^{\gamma^\ell_+} - \gamma^\ell_- (1-y^*)^{\gamma^\ell_-} }  \\
	\label{eqn:c2_def}
	c_2 &= \frac{\beta_\ell(1-y^*) - \tilde{h}(y^*) +\alpha_\ell\gamma^\ell_- (1-y^*)^{\gamma^\ell_-}}{\gamma^\ell_+ (1-y^*)^{\gamma^\ell_+} - \gamma^\ell_- (1-y^*)^{\gamma^\ell_-} } \,.
\end{align}
The boundary condition Equation \eqref{eqn:BCR} at $y = 1$ is automatically satisfied by the particular solution, and so the final coefficient $c_3$ is determined by the condition at $y = y^*$,
\begin{equation}
	\label{eqn:c3_def}
	c_3 = \frac{\beta_r(1-y^*)-\tilde{h}(y^*)}{\gamma^r_+ (1-y^*)^{\gamma^r_+}} \,.
\end{equation}
Finally, having determined the general solutions $V_\ell$ and $V_r$, the solution $V$ is obtained by enforcing continuity. This also fixes the value $y^*$, which must now satisfy the nonlinear transcendental equation
\begin{equation}
c_1 (1-y^*)^{\gamma^\ell_-}+c_2 (1-y^*)^{\gamma^\ell_+} + \alpha_\ell + \beta_\ell y^* = c_3 (1-y^*)^{\gamma^r_+} + \alpha_r + \beta_r y^* \,.
\end{equation}
By expanding the expression, it follows that
\begin{align}
	0=&\left((1-y^*)^{(\gamma^\ell_+ - \gamma^\ell_-)} - 1 \right)\left(\beta_\ell(1-y^*) - \tilde{h}(y^*)\right)
	-\alpha_\ell(\gamma^\ell_+ - \gamma^\ell_-) (1-y^*)^{\gamma^\ell_+}  \nonumber	\\
	&+\left(\gamma^\ell_+ (1-y^*)^{(\gamma^\ell_+ - \gamma^\ell_-)} - \gamma^\ell_-  \right)\left(\alpha_\ell-\alpha_r + (\beta_\ell-\beta_r )y^*- \frac{\beta_r(1-y^*)-\tilde{h}(y^*)}{\gamma^r_+} \right) \,.
\end{align}

\subsection{Comparison with $\mu \equiv 0$}

If the seller decides to take no action to enhance her reputation and long term wealth growth by letting $\mu \equiv 0$, then her value function solves 
\begin{equation}
\begin{aligned}
	-\frac{\sigma^2}{2}(1-y)^2\ V_{0}^{\prime \prime}(y) + \rho V_{0}(y) &= \tilde{h}(y) 	\qquad	 0<y<1 \\
	V_{0}(0) & = 0 \\
           V_{0}(1) & = \frac{1}{\rho} \tilde{h}(1) \,,
\end{aligned}
\end{equation}
which leads to a solution
\begin{equation}
	V_{0}(y) = \frac{\tilde{h}(y) }{\rho} - \frac{C}{\rho} (1-y)^{\frac{1}{2} + \sqrt{\frac{1}{4} + \frac{2 \rho}{\sigma^{2}} } } \,.
\end{equation}
It follows that  for $V$ which solves $(3.15)$ that 
\begin{equation}
\frac{A+C}{\rho} \epsilon =  V(1) - V_{0}(1) \leq \sup_{0 \leq y \leq 1} \Big | V(y) - V_{0}(y) \Big | \,.
\end{equation}
Hence, she enhances her overall expected wealth by order $\epsilon$ by pulsing strategies.

\section*{Acknowledgments}
The authors thank Andrew Christlieb, Nir Gavish, Song Yao, John Chadam, and Steven Shreve for useful discussions, as well as anonymous reviewers whose suggestions greatly improved this paper. Part of this work was done while Milan Bradonji\'c was at UCLA and LANL.

\bibliographystyle{acm}
\bibliography{reputation}
\end{document}